\documentclass[12pt]{amsart}
\usepackage{amssymb}
\usepackage{amscd}

\newcommand{\Q}{\mathbb{Q}}
\newcommand{\Z}{\mathbb{Z}}
\newcommand{\N}{\mathbb{N}}
\newcommand{\R}{\mathbb{R}}

\newcommand{\Ss}{\mathbb{S}}
\newcommand{\Aa}{\mathbb{A}}

\newcommand{\bbone}{\mathbf{1}}
\newcommand{\B}{\mathbf{B}}

\newcommand{\CalO}{\mathcal{O}}

\newcommand{\tH}{\widetilde{H}}
\newcommand{\WH}{W\! H}

\newcommand{\onemodd}{{1\, \mathrm{mod}\, d}}
\newcommand{\zeromodd}{{0\, \mathrm{mod}\, d}}
\newcommand{\nononemodd}{{\neq 1\, \mathrm{mod}\, d}}
\newcommand{\nonzeromodd}{{\neq 0\, \mathrm{mod}\, d}}

\newcommand{\Exp}{\mathrm{Exp}}

\newcommand{\Sech}{\mathrm{Sech}}
\newcommand{\sech}{\mathrm{sech}}
\newcommand{\Tanh}{\mathrm{Tanh}}
\newcommand{\Arcsinh}{\mathrm{Arcsinh}}
\newcommand{\arcsinh}{\mathrm{arcsinh}}
\newcommand{\Ind}{\mathrm{Ind}}

\newcommand{\tr}{\mathrm{tr}}

\newcommand{\rk}{\mathrm{rk}}
\newcommand{\ch}{\mathrm{ch}}
\newcommand{\sch}{\mathrm{sch}}
\newcommand{\sdim}{\mathrm{sdim}\,}

\newcommand{\isomto}{\overset{\sim}{\rightarrow}}

\newcommand{\noneq}{\natural}

\numberwithin{equation}{section}
\newtheorem{theorem}{Theorem}[section]
\newtheorem{corollary}[theorem]{Corollary}
\newtheorem{lemma}[theorem]{Lemma}
\newtheorem{proposition}[theorem]{Proposition}

\theoremstyle{definition}
\newtheorem{definition}[theorem]{Definition}

\newtheorem{example}[theorem]{Example}

\newcommand{\bdf}{\begin{definition}}
\newcommand{\edf}{\end{definition}\noindent}
\newcommand{\bex}{\begin{example}}
\newcommand{\eex}{\end{example}\noindent}
\newcommand{\bpr}{\begin{proposition}}
\newcommand{\epr}{\end{proposition}}
\newcommand{\blm}{\begin{lemma}}
\newcommand{\elm}{\end{lemma}}
\newcommand{\bth}{\begin{theorem}}
\renewcommand{\eth}{\end{theorem}}
\newcommand{\bpf}{\begin{proof}}
\newcommand{\epf}{\end{proof}\noindent}
\newcommand{\bcr}{\begin{corollary}}
\newcommand{\ecr}{\end{corollary}\noindent}
\newcommand{\beq}{\begin{equation}}
\newcommand{\eeq}{\end{equation}}
\newcommand{\bes}{\begin{equation*}}
\newcommand{\ees}{\end{equation*}}
\newcommand{\ben}{\begin{enumerate}}
\newcommand{\een}{\end{enumerate}}
\begin{document}
\title[Homology of sub-posets of Dowling lattices]
{Plethysm for wreath products and homology of sub-posets of Dowling lattices}
\author{Anthony Henderson}
\address{School of Mathematics and Statistics,
University of Sydney, NSW 2006, AUSTRALIA}
\email{anthonyh@maths.usyd.edu.au}
\thanks{This work was supported by Australian Research Council grant DP0344185}
\begin{abstract}
We prove analogues for sub-posets of the Dowling lattices of the
results of Calderbank, Hanlon, and Robinson on homology of sub-posets
of the partition lattices. The technical tool used is the wreath product
analogue of the tensor species of Joyal.
\end{abstract}
\maketitle
\section*{Introduction}
For any positive integer $n$ and finite group $G$, the Dowling lattice
$Q_n(G)$ is a poset with an action of the wreath product group
$G\wr S_n$. If $G$ is trivial, $Q_n(\{1\})$ can be identified with the
partition lattice $\Pi_{n+1}$ (on which $S_n$ acts as a subgroup of $S_{n+1}$).
If $G$ is the cyclic group of order $r$ for $r\geq 2$, $Q_n(G)$ can be
identified with the lattice of intersections of reflecting hyperplanes in the
reflection representation of $G\wr S_n$. For general $G$, the underlying
set of $Q_n(G)$ can be thought of as
the set of all pairs $(I,\pi)$ where $I\subseteq\{1,\cdots,n\}$ and
$\pi$ is a set partition of $G\times(\{1,\cdots,n\}\setminus I)$ whose
parts $G$ permutes freely; see Definition \ref{dowlingdef} below for
the partial order.

In Section 1 we will 
define various sub-posets $P$ of $Q_n(G)$, containing the
minimum element $\hat{0}$ and the maximum element $\hat{1}$, 
which are stable under the
action of $G\wr S_n$. For completeness' sake we include the cases of $Q_n(G)$
itself and two other sub-posets which have been studied before,
but the main interest lies in two new families of sub-posets, defined
using a fixed integer $d\geq 2$:
$Q_n^{\onemodd}(G)$, given by the congruence conditions
$|I|\equiv 0$ mod $d$ and $|K|\equiv 1$ mod $d$ for all parts $K$ of $\pi$,
and $Q_n^{\zeromodd}(G)$, given by the condition
$|K|\equiv 0$ mod $d$ for all parts $K$ of $\pi$. These definitions are
modelled on those of the sub-posets $\Pi_n^{(1,d)}$ and $\Pi_n^{(0,d)}$
of the partition lattice studied by 
Calderbank, Hanlon, and Robinson in \cite{chr}.
We will prove that all our
sub-posets $P$ are pure (i.e.\ graded) and Cohen-Macaulay, so
the only non-vanishing reduced homology group of 
$P\setminus\{\hat{0},\hat{1}\}$ is the top homology
$\tH_{l(P)-2}(P\setminus\{\hat{0},\hat{1}\};\Q)$.
We take rational coefficients so that we can regard this homology
as a representation of $G\wr S_n$ over $\Q$.

The main aim of this paper is to find in each case a formula for the character of
this representation, analogous to the formulae proved in \cite{chr}. 
The last paragraph of that paper
hoped specifically for a Dowling lattice analogue of \cite[Theorem 6.5]{chr},
but our Theorem \ref{0moddthm} casts doubt on its existence. 
(In the cases of the
previously-studied posets, we recover Hanlon's formula from \cite{hanlonwreath}
and other results which were more or less known.)

In \cite{rains}, Rains has
applied \cite[Theorem 4.7]{chr} to compute
the character of $S_n$ on the cohomology of the manifold
$\overline{\mathcal{M}_{0,n}}(\R)$ (the real points of the moduli space
of stable genus $0$ curves with $n$ marked points). His suggestion
that a Dowling lattice analogue would have a similar application
to the cohomology of the real points of De Concini-Procesi compactifications
was the motivation for this work. For more detail on the connection,
see the remarks following \eqref{noneqmainbneqn}.

In Section 2 we recall the combinatorial framework used by
Macdonald to write down characters of representations of wreath
products, and state our results. In Section 3 we introduce
the functorial concept of a $(G\wr\Ss)$-module, a generalization
of Joyal's notion of tensor species; this concept comes from
\cite{mywreath}, and we recall the connection proved there
with generalizations of plethysm. In Section 4 we use this
technology, and the `Whitney homology method' of Sundaram, to prove our results. In Section 5 we extend the results to the setting of
Whitney homology, thus computing the `equivariant characteristic polynomials'
of our posets.
\section{Some Cohen-Macaulay sub-posets of Dowling lattices}
In this section we define the Dowling lattices and the sub-posets
of interest to us, and prove that they are Cohen-Macaulay.
A convenient reference
for the basic definitions and techniques of Cohen-Macaulay posets 
is \cite{wachs}; the key result for us is the Bj\"orner-Wachs criterion,
\cite[Theorem 4.2.2]{wachs} (proved in \cite{bjornerwachs}),
that a pure bounded poset with a recursive atom ordering is Cohen-Macaulay.
For any nonnegative integer $n$, write $[n]$ for $\{1,\cdots,n\}$
(so $[0]$ is the empty set), and $S_n$ for the symmetric group
of permutations of $[n]$. 

For any finite set $I$,
let $\Pi(I)$ denote the poset of partitions of the set $I$,
where a \textbf{partition} $\pi$ of $I$ is a set of nonempty disjoint subsets
of $I$ whose union is $I$. These subsets $K\in\pi$ are referred to as
the \textbf{parts} of $\pi$. The partial order on $\Pi(I)$ is by refinement;
$\Pi(I)$ is a geometric lattice, isomorphic to $\Pi([|I|])$ which is more
commonly written $\Pi_{|I|}$. (We use the convention that the
empty set has a single partition, which as a set is itself empty. Therefore
$\Pi(\emptyset)=\Pi_0$ is a one-element poset, like $\Pi_1$.) 

Fix a finite group $G$, and view the wreath product $G\wr S_n$ as
the group of permutations of $G\times [n]$ which commute with the
action of $G$ (by left multiplication on the first factor). Our definition
of the corresponding Dowling lattice is as follows.
\bdf \label{dowlingdef}
For $n\geq 1$, let $Q_n(G)$ be the poset
of pairs $(J,\pi)$ where $J$ is a $G$-stable subset of $G\times [n]$ and
$\pi\in\Pi((G\times [n])\setminus J)$ is such that 
$G$ permutes its parts freely, i.e.\ for all $1\neq g\in G$ and
$K\in\pi$, $K\neq g.K\in\pi$.
The partial order on these pairs is defined so that $(J,\pi)\leq(J',\pi')$
is equivalent to the following two conditions:
\ben
\item $J\subseteq J'$, and
\item for all parts $K\in\pi$, either $K\subseteq J'$ or $K$ is contained
in a single part of $\pi'$.
\een
\edf
We have an obvious action of $G\wr S_n$ on the poset $Q_n(G)$.
Of course, $J$ must be of the form $G\times I$ for some subset 
$I\subseteq [n]$, so we could just as well have used $I$ in the definition,
as in the introduction;
once one has taken into account this and other such variations,
it should be clear that $Q_n(G)$ is isomorphic, as a $(G\wr S_n)$-poset,
to Dowling's original lattice in \cite{dowling}
and to the various alternative definitions given
in \cite{hanlonwreath}, \cite{gottliebwachs}, and \cite{hultman}.
(The justification for adding yet another definition to the list
will come when we adopt a functorial point of view.)
The minimum element $\hat{0}$
is the pair $(\emptyset,\{\{(g,m)\}\,|\,g\in G, m\in[n]\})$, 
and the maximum element
$\hat{1}$ is the pair $(G\times [n],\emptyset)$.
Dowling proved in \cite{dowling} that $Q_n(G)$ is a geometric
lattice, and hence it is Cohen-Macaulay; its rank function is
\beq \label{dowlingrankeqn}
\rk(J,\pi)=n-\frac{|\pi|}{|G|},
\eeq
so the length of the lattice as a whole is $n$.

Special cases of this lattice are more familiar.
Clearly $Q_n(\{1\})\cong\Pi(\{0,1,\cdots,n\})$ via the map
which sends $(J,\pi)$ to $\{J\cup\{0\}\}\cup\pi$;
and $Q_n(\{\pm 1\})$ is the signed partition lattice, also known as
the poset of (conjugate) parabolic subsystems of a root system of type $B_n$.
More generally, when $G$ is cyclic of order $r\geq 2$,
$Q_n(G)$ can be
identified with the lattice of intersections of reflecting hyperplanes in the
reflection representation of $G\wr S_n$, i.e.\ the lattice denoted
$L(\mathcal{A}_n(r))$ in \cite[\S 6.4]{orlikterao}.

Before we define the sub-posets we are mainly interested in,
let us also consider two sub-posets given by a condition on $J$:
\bdf
For $n\geq 1$, let $R_n(G)$ be the sub-poset
of $Q_n(G)$ consisting of pairs $(J,\pi)$ where either $J=\emptyset$ or 
$J=G\times [n]$. For $n\geq 2$ and assuming that $G\neq\{1\}$, 
let $Q_n^\sim(G)$ be the sub-poset
of $Q_n(G)$ consisting of pairs $(J,\pi)$ where $\frac{|J|}{|G|}\neq 1$.
\edf
Clearly the minimum and maximum elements of $Q_n(G)$ are in $R_n(G)$
(indeed, we allow $J=G\times [n]$ merely in order to include $\hat{1}$);
likewise for $Q_n^\sim(G)$, given that $n\geq 2$.
Since $R_n(G)\setminus\{\hat{1}\}$ is a lower order ideal of $Q_n(G)$,
it is obvious that $R_n(G)$ is pure of length $n$, with rank function
again given by \eqref{dowlingrankeqn}; a moment's thought shows that
the same is true for $Q_n^\sim(G)$.
It is also easy to see that $R_n(G)\setminus\{\hat{1}\}$ is a geometric
semilattice in the sense of Wachs and Walker 
(see \cite[Definition 4.2.6]{wachs}), so $R_n(G)$ is Cohen-Macaulay 
by \cite[Theorem 4.2.7]{wachs}. An alternative proof of this is
provided by \cite[Corollary 3.12]{hultman}, where $R_n(G)\setminus\{\hat{0},
\hat{1}\}$ is called $\underline{\Pi}_n^G$. Note that
$R_n(\{1\})\setminus\{\hat{1}\}\cong\Pi_n$, so $R_n(\{1\})$ is
$\Pi_n$ with an extra maximum element adjoined. 
One can also interpret $R_n(\{\pm 1\})\setminus\{\hat{1}\}$
as the poset of (conjugate) parabolic subsystems of a root system
of type $B_n$ all of whose components are of type $A$.
As for $Q_n^\sim(G)$, note that it is closed under the join operation
of $Q_n(G)$, and two elements of $Q_n^\sim(G)$ have a meet in
$Q_n^\sim(G)$ which is $\leq$ their meet in $Q_n(G)$. The assumption
that $G\neq\{1\}$ ensures that every element of $Q_n^\sim(G)$
is a join of atoms, so
$Q_n^\sim(G)$ is another geometric lattice, and hence it is Cohen-Macaulay.
When $G$ is cyclic of order $r\geq 2$, $Q_n^\sim(G)$ can be
identified with the lattice denoted
$L(\mathcal{A}_n^0(r))$ in \cite[Section 6.4]{orlikterao}; for instance,
$Q_n^\sim(\{\pm 1\})$ is
the poset of (conjugate) parabolic subsystems of a root system of type $D_n$.

Now we turn to the analogues of the sub-posets of the partition lattices
considered by Calderbank, Hanlon, and Robinson.
\bdf \label{1modddef}
For $n\geq 1$ and $d\geq 2$, let $Q_n^\onemodd(G)$ be the sub-poset
of $Q_n(G)$ consisting of pairs $(J,\pi)$ satisfying the following conditions:
\ben
\item for all $K\in\pi$, $|K|\equiv 1$ mod $d$; and
\item either $\frac{|J|}{|G|}\equiv 0$ mod $d$ or $J=G\times [n]$.
\een
\edf
Note that when $n\equiv 0$ mod $d$, there is no need to
explicitly allow $J=G\times [n]$; otherwise, 
allowing this has the effect of ensuring
that the maximum element $\hat{1}$ is included. Clearly the minimum element
$\hat{0}$ of $Q_n(G)$ also belongs to $Q_n^\onemodd(G)$. To explain
the congruence condition in (2), note that under the isomorphism
$Q_n(\{1\})\cong\Pi(\{0,1,\cdots,n\})\cong\Pi_{n+1}$, 
$Q_n^\onemodd(\{1\})$ corresponds to the `$1$ mod $d$ partition lattice'
$\Pi_{n+1}^{(1,d)}$ considered in \cite{chr}. Also, one can interpret
$Q_n^\onemodd(\{\pm 1\})\setminus\{\hat{1}\}$ as the poset of proper
(conjugate) parabolic subsystems of a root system
of type $B_n$ all of whose components have rank divisible by $d$; this
poset has arisen in recent work of Rains, who has proved the
following result in that case.
\bpr
$Q_n^\onemodd(G)$ 
is a totally semimodular pure poset with rank function
\[ \rk(J,\pi)=\left\{\begin{array}{cl}
\lceil\frac{n}{d}\rceil,&\text{ if $J=G\times n$,}\\
\frac{n}{d}-\frac{|\pi|}{d|G|},&\text{ otherwise.}
\end{array}\right. \]
Its length is $\lceil\frac{n}{d}\rceil$.
\epr
\bpf
If $(J',\pi')$ covers $(J,\pi)$ in $Q_n^\onemodd(G)$, then there are
two possibilities:
\ben
\item $J'=J$, in which case $\pi'$ must be obtained from $\pi$
by merging $(d+1)$ $G$-orbits of parts into a single $G$-orbit
of parts, \textbf{or}
\item $J'\supset J$, in which case $J'$ must be the union of $J$
together with $d$ $G$-orbits of parts of $\pi$ (or 
$n-d\lfloor\frac{n}{d}\rfloor$ $G$-orbits, if $J'=G\times [n]$ and 
$n\not\equiv 0$ mod $d$).
\een
In either case one sees immediately that the purported rank
of $(J',\pi')$ is one more than that of $(J,\pi)$, so
this is indeed the rank function, and $Q_n^\onemodd(G)$ is pure.
To show that $Q_n^\onemodd(G)$ is totally semimodular (see \cite[4.2]{wachs}),
it suffices to check the condition at $\hat{0}$, since for every
$(J,\pi)\in Q_n^\onemodd(G)$, the principal upper order ideal
$[(J,\pi),\hat{1}]$ is isomorphic to $Q_{|\pi|/|G|}^\onemodd(G)$.
That is, we need only prove the following: if $a$ and $b$ are distinct atoms
of $Q_n^\onemodd(G)$, $a\vee b$ is their join in the lattice $Q_n(G)$,
and $c\in Q_n^\onemodd(G)$ satisfies $c\geq a\vee b$ and is minimal
with this property in $Q_n^\onemodd(G)$, then $\rk(c)=2$ for the rank
function we have just found. Since there are two types of atoms corresponding
to the two kinds of covering relation, we have several 
cases to consider.\newline
\textbf{Case 1:} $a=(\emptyset,\pi_a)$ and $b=(\emptyset,\pi_b)$. 
Let $A$ be the union of the non-singleton parts
of $\pi_a$, which all have size $d+1$ and form a single $G$-orbit. 
Define $B$ similarly for $\pi_b$.
Let $\pi_a\vee\pi_b$ denote the join in $\Pi(G\times [n])$.
\newline 
\textbf{Subcase 1a:} $\frac{|A\cap B|}{|G|}=0$ or $1$. 
Then $a\vee b=(\emptyset,\pi_a\vee\pi_b)$ is itself in $Q_n^\onemodd(G)$, so
$c=a\vee b$ and $\rk(c)=2$.\newline
\textbf{Subcase 1b:} $\frac{|A\cap B|}{|G|}\geq 2$, and
$a\vee b=(\emptyset,\pi_a\vee\pi_b)$.
(This means that the mergings of $A$ and of $B$
are `compatible' on the overlap.) Then $\pi_a\vee\pi_b$ has a unique
$G$-orbit of non-singleton parts, whose union is $A\cup B$.
Since $A\neq B$, we have $d+2\leq \frac{|A\cup B|}{|G|}
=2d+2-\frac{|A\cap B|}{|G|}\leq 2d$.
Then either $c=(J,\hat{0}_{\Pi((G\times [n])\setminus J)})$ where
$J\supseteq A\cup B$ has size $\min\{2d|G|,n|G|\}$, or
$c=(\emptyset,\pi_c)$ where $\pi_c$ has a unique $G$-orbit
of non-singleton parts, all of size $2d+1$, whose union contains $A\cup B$.
In either case $\rk(c)=2$.\newline
\textbf{Subcase 1c:} $\frac{|A\cap B|}{|G|}\geq 2$, 
and $a\vee b=(A\cup B,\hat{0}_{\Pi((G\times [n])\setminus(A\cup B))})$.
(This means that the mergings of $A$ and of $B$ are `not compatible' on the
overlap, as can happen when $G$ is non-trivial). We have
$d+1\leq\frac{|A\cup B|}{|G|}
=2d+2-\frac{|A\cap B|}{|G|}\leq 2d$,
so $c$ must be of the form $(J,\hat{0}_{\Pi((G\times [n])\setminus J)})$ where
$J\supseteq A\cup B$ has size $\min\{2d|G|,n|G|\}$. Thus $\rk(c)=2$.\newline
\textbf{Case 2:} the atoms $a$ and $b$ are of different types.
Without loss of generality, assume 
$a=(A,\hat{0}_{\Pi((G\times [n])\setminus A)})$ where $|A|=d|G|$,
and $b=(\emptyset,\pi_b)$ for $B$ as above.\newline
\textbf{Subcase 2a:} $A\cap B=\emptyset$. 
Then $a\vee b=(A,\pi_b|_{(G\times [n])\setminus A})$ is itself in
$Q_n^\onemodd(G)$, so $c=a\vee b$ and $\rk(c)=2$.\newline
\textbf{Subcase 2b:} $A\cap B\neq\emptyset$. Then 
$a\vee b=(A\cup B,\hat{0}_{\Pi((G\times [n])\setminus(A\cup B))})$, and
$d+1\leq\frac{|A\cup B|}{|G|}
=2d+1-\frac{|A\cap B|}{|G|}\leq
2d$, so $c$ must be as in Subcase 1c.\newline
\textbf{Case 3:} $a=(A,\hat{0}_{\Pi((G\times [n])\setminus A)})$,
$b=(B,\hat{0}_{\Pi((G\times [n])\setminus B)})$ where
$|A|=|B|=d|G|$. Then 
$a\vee b=(A\cup B,\hat{0}_{\Pi((G\times [n])\setminus(A\cup B))})$.
Since $A\neq B$, we have
$d+1\leq\frac{|A\cup B|}{|G|}
=2d-\frac{|A\cap B|}{|G|}\leq
2d$, so $c$ must be as in Subcase 1c.
\epf
We deduce via \cite[Theorem 4.2.3]{wachs} that $Q_n^\onemodd(G)$ is 
Cohen-Macaulay.

Finally, we consider the analogue of the `$d$-divisible partition lattice'.
\bdf
For $n\geq 1$ and $d\geq 2$, 
let $Q_n^\zeromodd(G)$ be the sub-poset
of $Q_n(G)$ consisting of pairs $(J,\pi)$ such that either
\ben
\item $(J,\pi)$ is the minimum element of $Q_n(G)$, i.e. $J=\emptyset$ and
$|K|=1$ for all $K\in\pi$, or
\item $|K|\equiv 0$ mod $d$ for all $K\in\pi$.
\een
\edf
Note that the maximum element of $Q_n(G)$ vacuously satisfies condition (2),
so this poset is certainly bounded. If $n\equiv -1$ mod $d$,
then under the isomorphism
$Q_n(\{1\})\cong\Pi(\{0,1,\cdots,n\})\cong\Pi_{n+1}$, 
$Q_n^\zeromodd(\{1\})$ corresponds to the poset
$\Pi_{n+1}^{(0,d)}$ considered in \cite{chr}. If $n\not\equiv -1$ mod $d$,
then $Q_n^\zeromodd(\{1\})$ does not correspond to anything in \cite{chr}.
\bpr
$Q_n^\zeromodd(G)$ is a pure join semilattice and has a 
recursive atom ordering. Its rank
function is
\[ \rk(J,\pi)=\left\{\begin{array}{cl}
0,&\text{ if $(J,\pi)=\hat{0}$,}\\
\lfloor\frac{n}{d}\rfloor+1-\frac{|\pi|}{|G|},&\text{ otherwise.}
\end{array}\right. \]
Its length is $\lfloor\frac{n}{d}\rfloor+1$.
\epr
\bpf
It is obvious that $Q_n^\zeromodd(G)\setminus\{\hat{0}\}$ is an upper
order ideal of $Q_n(G)$, so
$Q_n^\zeromodd(G)$ is a join semilattice.
The atoms of $Q_n^\zeromodd(G)$ are those $(J,\pi)$
where $|K|=d$ for all $K\in\pi$ and 
$\frac{|\pi|}{|G|}=\lfloor\frac{n}{d}\rfloor$; these all have rank
$(n-\lfloor\frac{n}{d}\rfloor)$ as elements of $Q_n(G)$, so
$Q_n^\zeromodd(G)$ is pure and has the claimed rank function.
To find a recursive atom ordering
(see \cite[Definition 4.2.1]{wachs}), note that for any non-mimimum element
$(J,\pi)$ of $Q_n^\zeromodd(G)$, the principal upper order ideal
$[(J,\pi),\hat{1}]$ is totally semimodular, being
isomorphic to $Q_{|\pi|/|G|}(G)$. Thus we need only check
that the atoms of $Q_n^\zeromodd(G)$
can be ordered $a_1,\cdots,a_t$ so that: 
\beq \label{orderingcondeqn}
a_i,a_j<y, i<j\Rightarrow
\exists z\leq y, z\text{ covers }a_j\text{ and }
a_k, \exists k<j.
\eeq
Such an ordering (inspired by \cite[Exercise 4.3.6(a)]{wachs})
can be defined as follows. For each $d$-element
subset $I\subset[n]$, let $\Psi(I)$ be the set of
partitions of $G\times I$ on whose parts $G$ acts freely and transitively
(there are $|G|^{d-1}$ elements in this set). An atom
$(J,\pi)$ of $Q_n^\zeromodd(G)$ is uniquely determined by the following data:
\ben
\item a $(n-d\lfloor\frac{n}{d}\rfloor)$-element subset $I_0$ of $[n]$, such
that $J=G\times I_0$; and
\item a partition of $[n]\setminus I_0$ into $d$-element subsets
$I_1,\cdots,I_{\lfloor\frac{n}{d}\rfloor}$, each $I_s$ equipped with
a partition $\psi_s\in\Psi(I_s)$, such that $\pi=\bigcup_{s}\psi_s$.
\een
From these data, construct a word by concatenating the elements of $I_0$
(in increasing order) followed by the elements of $I_1$
(in increasing order), $I_2$
(in increasing order), and so on up to $I_{\lfloor\frac{n}{d}\rfloor}$, where
the ordering of $I_1,\cdots,I_{\lfloor\frac{n}{d}\rfloor}$ themselves
is determined by the order of their smallest elements. Then order the
atoms by lexicographic order of these words; within atoms with the same word,
use the order given by some arbitrarily chosen orderings of the sets
$\Psi(I)$ for all $d$-element subsets $I$ (applied lexicographically, so
the ordering of $\Psi(I_1)$ is applied first, then in case of equality of
$\psi_1$ the ordering of $\Psi(I_2)$ is applied, etc.).

We now prove that this ordering satisfies the condition 
\eqref{orderingcondeqn}. Let $a_j=(J,\pi)$ have associated
$I_s$ and $\psi_s$ as above, let $y=(J',\pi')\in Q_n^\zeromodd(G)$ be
such that $(J',\pi')>(J,\pi)$, and suppose that $(J',\pi')$ is not greater than
any common cover of $a_j$ and an earlier atom. We must deduce from this
that $(J,\pi)$ is the earliest atom which is $<(J',\pi')$. Firstly,
let $K$ be any part of $\pi'$, and let $s_1<\cdots<s_t$ be such that
$\bigcup_{g\in G}g.K=G\times (I_{s_1}\cup\cdots\cup I_{s_t})$.
Suppose that for some $i$, $a=\max(I_{s_i})>\min(I_{s_{i+1}})=b$. Let
$g_a,g_b\in G$ be such that $(g_a,a),(g_b,b)\in K$. There is an
element $w\in G\wr S_n$ defined by
\[ w.(g,c)=\left\{\begin{array}{cl}
(gg_a^{-1}g_b,b),&\text{ if $c=a$,}\\
(gg_b^{-1}g_a,a),&\text{ if $c=b$,}\\
(g,c),&\text{ if $c\neq a,b$.}
\end{array}\right. \]
It is clear that $w.(J,\pi)$ is an earlier atom than $(J,\pi)$, and their
join is a common cover which is $\leq (J',\pi')$, contrary to assumption.
Hence we must have $\max(I_{s_i})<\min(I_{s_{i+1}})$ for all $i$.
Thus the parts of $\pi$ contained in $K$ are simply those which one obtains
by ordering the elements of $K$ by their second component, and chopping
that list into $d$-element sublists. By similar arguments (details omitted),
one can show that $I_0$ must consist of the $(n-d\lfloor\frac{n}{d}\rfloor)$
smallest numbers occurring in the second components of elements of $J'$,
and that the parts of $\pi$ contained in $J'$ are those obtained by
listing the remaining such numbers in increasing order, chopping that list
into $d$-element sublists, and choosing for each resulting $I_s$ the
smallest element of $\Psi(I_s)$ (for the fixed order on this set). It is
clear from this construction of $(J,\pi)$ that it is 
the earliest atom which is $<(J',\pi')$.
\epf
We deduce via \cite[Theorem 4.2.2]{wachs} that $Q_n^\zeromodd(G)$ is 
Cohen-Macaulay.
\section{Statement of the results}
In this section, after introducing some necessary notation,
we state our results on the character of $G\wr S_n$ on
$\tH_{l(P)-2}(\overline{P};\Q)$
for each of the sub-posets $P$ of $Q_n(G)$ defined in the previous section;
here $\overline{P}$ denotes the `proper part' $P\setminus\{\hat{0},\hat{1}\}$.
Since the posets are Cohen-Macaulay, this is the only reduced homology
group of $\overline{P}$ which can be nonzero.
Hence $\dim\tH_{l(P)-2}(\overline{P};\Q)
=(-1)^{l(P)}\mu(P)$.
(We follow the usual convention that $\tH_{-1}(\emptyset;\Q)$ is 
one-dimensional.)

Let $G_*$ denote the set of conjugacy classes of $G$. Following
\cite[Chapter I, Appendix B]{macdonald}, we introduce
the polynomial ring $\Lambda_G:=\Q[p_i(c)]$
in indeterminates $p_i(c)$, one for each positive integer $i$ and
conjugacy class $c\in G_*$. This ring is $\N$-graded by setting
$\deg(p_i(c))=i$.
The character of a representation $M$ of $G\wr S_n$ over $\Q$
is encapsulated in its \textbf{Frobenius characteristic}
\beq
\ch_{G\wr S_n}(M):=\frac{1}{|G|^n n!}\sum_{x\in G\wr S_n}
\tr(x,M)\Psi(x),
\eeq
which is a homogeneous element of $\Lambda_G$ of degree $n$.
The definition of the \textbf{cycle index} $\Psi(x)$
is $\prod_{i\geq 1,c\in G_*}p_i(c)^{a_i(c)}$ if $x$ lies in the
conjugacy class of elements with $a_i(c)$ cycles of length $i$ and type $c$
(see [loc.\ cit.]). Recall that such an element has centralizer of order
\beq
\prod_{\substack{i\geq 1\\c\in G_*}} \left(\frac{|G|}{|c|}i\right)^{a_i(c)} 
a_i(c)!,
\eeq
as shown in [loc.\ cit., (3.1)]; its trace on $M$ can be recovered
from $\ch_{G\wr S_n}(M)$ by multiplying the coefficient of
$\prod_{i\geq 1,c\in G_*}p_{i(c)}^{a_i(c)}$ by this centralizer order.
For any $f\in\Lambda_G$,
we write $f^\noneq$ for its `non-equivariant specialization',
the element of $\Q[x]$ obtained from $f$ by setting $p_1(\{1\})$ to $x$
and all other $p_i(c)$ to $0$. Clearly
\beq
\ch_{G\wr S_n}(M)^\noneq = (\dim M)\frac{x^n}{|G|^n n!}.
\eeq
When $G=\{1\}$ we write $p_i$ for $p_i(\{1\})$, as usual in the theory
of symmetric groups and symmetric functions.

Now it is well known that $\Lambda_{\{1\}}$ has an associative operation
called \textbf{plethysm}, for which $p_1$ is an identity. Less well known
is that $\Lambda_{G}$ has a pair of `plethystic actions' of $\Lambda_{\{1\}}$,
one on the left and one on the right; in the terminology of 
\cite{borgerwieland}, $\Lambda_{\{1\}}$
is a plethory, and $\Lambda_G$ is a 
$\Lambda_{\{1\}}$--$\Lambda_{\{1\}}$--biring.
The left plethystic action is an operation
$\circ:\Lambda_{\{1\}}\times\Lambda_{G}\to\Lambda_G$, which is 
uniquely defined by:
\begin{enumerate}
\item for all $g\in \Lambda_{G}$, the map 
$\Lambda_{\{1\}}\to \Lambda_G:f\mapsto f\circ g$
is a homomorphism of $\Q$-algebras;
\item for any $i\geq 1$, the map
$\Lambda_{G}\to\Lambda_G:g\mapsto p_{i}\circ g$ is a 
homomorphism of $\Q$-algebras;
\item $p_i\circ p_j(c)=p_{ij}(c)$.
\end{enumerate}
This action implicitly appears in \cite{polyfunctors}. 
The more interesting right plethystic action, made explicit
for the first time in \cite[Section 5]{mywreath}, is an operation
$\circ:\Lambda_G\times\Lambda_{\{1\}}\to\Lambda_G$, 
which is uniquely defined by:
\begin{enumerate}
\item for all $g\in\Lambda_{\{1\}}$, the map 
$\Lambda_G\to\Lambda_G:f\mapsto f\circ g$
is a homomorphism of $\Q$-algebras;
\item for any $i\geq 1$, $c\in G_*$, the map
$\Lambda_{\{1\}}\to\Lambda_G:g\mapsto p_{i}(c)\circ g$ is a 
homomorphism of $\Q$-algebras;
\item $p_{i}(c)\circ p_j=p_{ij}(c^j)$, where $c^j$ denotes the conjugacy class
of $j$th powers of elements of $c$.
\end{enumerate}
If $G=\{1\}$ both these actions become the usual operation of plethysm.
We have $(f\circ g)\circ h=f\circ(g\circ h)$ whenever $f,g,h$ live in
the right combination of $\Lambda_{\{1\}}$ and/or $\Lambda_G$
for both sides to be defined; moreover, $p_1\circ f=f\circ p_1=f$ for all
$f\in\Lambda_G$. Note that under the non-equivariant
specialization, all cases of $\circ$ become simply the
substitution of one polynomial in $\Q[x]$ into another. For more on the
`meaning' of these plethystic actions, consult \cite[Section 5]{mywreath}
or the next section.

Since our formulae use generating functions which combine
$(G\wr S_n)$-modules for infinitely many $n$,
we need to enlarge $\Lambda_G$ to the formal power series ring 
$\Aa_G:=\Q[\![p_i(c)]\!]$, which we give its usual topology (coming from the
$\N$-filtration). We extend the non-equivariant specialization
in the obvious way (that is, by continuity): for 
$f\in\Aa_G$, $f^\noneq$ is an element
of the formal power series ring $\Q[\![x]\!]$. Just as one cannot substitute
a formal power series with nonzero constant term into another formal power
series, the extensions of $\circ$ to this context require a slight
restriction. Let $\Aa_{G,+}$ be the ideal of $\Aa_G$ consisting 
of elements whose degree-$0$ term vanishes. Then the left plethystic 
action extends to an operation
$\circ:\Aa_{\{1\}}\times\Aa_{G,+}\to\Aa_G$, and the right plethystic
action extends to an operation 
$\circ:\Aa_{G}\times\Aa_{\{1\},+}\to\Aa_G$. Of course, the associativity and
identity properties continue to hold.

An important element of $\Aa_G$ is the sum of the characteristics of the
trivial representations:
\bes
\Exp_G:=\sum_{n\geq 0}\ch_{G\wr S_n}(\bbone)
=\sum_{(a_i(c))}\prod_{\substack{i\geq 1\\c\in G_*}}
\frac{(\frac{|c|}{|G| i}p_i(c))^{a_i(c)}}{a_i(c)!}
=\exp(\sum_{\substack{i\geq 1\\c\in G_*}}\frac{|c| p_i(c)}{|G| i}).
\ees
Clearly $\Exp_G^\noneq=\exp(\frac{x}{|G|})$. We write $\Exp_{\{1\}}=
\exp(\sum_{i\geq 1}\frac{p_i}{i})$ simply as $\Exp$; one has
\beq
\Exp_G=\Exp\circ(\sum_{c\in G_*}\frac{|c|}{|G|}p_1(c)).
\eeq
It is well known that
the plethystic inverse of $\Exp-1$ in $\Aa_{\{1\},+}$ is
\beq
L:=\sum_{d\geq 1}\frac{\mu(d)}{d}\log(1+p_d).
\eeq
In other words, $L\circ(\Exp-1)=(\Exp-1)\circ L=p_1$.
With this notation, the famous result of Stanley that
$\tH_{n-3}(\overline{\Pi_n};\Q)\cong\varepsilon_n\otimes
\Ind_{\mu_n}^{S_n}(\psi)$, where
$\psi$ is a faithful character of the cyclic group $\mu_n$, 
can be rephrased as
\beq \label{stanleyeqn}
p_1+\sum_{n\geq 2}(-1)^{n-1}\ch_{S_n}(\tH_{n-3}(\overline{\Pi_n};\Q))=L.
\eeq
(The proof of this fact will be recalled in Section 4.)

We can rephrase \cite[Corollary 2.2]{hanlonwreath}
in a similar way.
\bth \label{hanlonthm} \textup{(Hanlon)} In $\Aa_G$ we have the equation
\[ 1+\sum_{n\geq 1}(-1)^n \ch_{G\wr S_n}(\tH_{n-2}
(\overline{Q_n(G)};\Q))=(\Exp_G\circ L)^{-1}. \]
\eth
\noindent
We will give a new proof of Theorem \ref{hanlonthm} 
in Section 4. To see that this is
equivalent to Hanlon's statement, note that
\bes
\begin{split}
(\Exp_G\circ L)^{-1}&=
\exp(-\sum_{\substack{i\geq 1\\c\in G_*}}\frac{|c| p_i(c)}{|G| i})
\circ(\sum_{d\geq 1}\frac{\mu(d)}{d}\log(1+p_d))\\
&=\exp(-\sum_{\substack{i\geq 1\\d\geq 1\\c\in G_*}}
\frac{|c|\mu(d)}{|G|id}\log(1+p_{id}(c^d)))\\
&=\prod_{\substack{l\geq 1\\c\in G_*}} (1+p_l(c))^{F(l,c)},
\end{split}
\ees
where
\beq
F(l,c):=-\frac{1}{|G|l}\sum_{d\mid l}\mu(d)|\{g\in G\,|\,g^d\in c\}|,
\eeq
which is easily equated with Hanlon's $F(l,c,1)$.
Applying ${}^\noneq$ to Theorem \ref{hanlonthm}, we derive the 
non-equivariant version:
\beq
1+\sum_{n\geq 1}(-1)^n \dim \tH_{n-2}
(\overline{Q_n(G)};\Q)\frac{x^n}{|G|^n n!}=(1+x)^{-1/|G|},
\eeq
which is equivalent to the well-known fact that
\beq
\dim \tH_{n-2}(\overline{Q_n(G)};\Q)
=(|G|+1)(2|G|+1)\cdots((n-1)|G|+1).
\eeq
Note that the $G=\{1\}$ special case of Theorem \ref{hanlonthm} is
\bes
1+\sum_{n\geq 1}(-1)^n \ch_{S_n}(\tH_{n-2}
(\overline{Q_n(\{1\})};\Q))=(1+p_1)^{-1},
\ees
which can also be obtained by applying $\frac{\partial}{\partial p_1}$ 
to both sides of \eqref{stanleyeqn}.

For the `all or nothing' sub-poset $R_n(G)$, we have the following result,
to be proved in Section 4.
\bth \label{secondthm}
In $\Aa_G$ we have the equation
\[ \sum_{n\geq 1}(-1)^n \ch_{G\wr S_n}(\tH_{n-2}
(\overline{R_n(G)};\Q))=1-\Exp_G\circ L. \]
\eth
\noindent
The non-equivariant version is
\beq
\sum_{n\geq 1}(-1)^n \dim \tH_{n-2}
(\overline{R_n(G)};\Q)\frac{x^n}{|G|^n n!}=1-(1+x)^{1/|G|},
\eeq
which is equivalent to the result of Hultman (\cite[Corollary 3.12]{hultman}):
\beq
\dim \tH_{n-2}
(\overline{R_n(G)};\Q)=(|G|-1)(2|G|-1)\cdots((n-1)|G|-1).
\eeq
Note also that the $G=\{1\}$ case of Theorem \ref{secondthm} is
\bes
\sum_{n\geq 1}(-1)^n \ch_{S_n}(\tH_{n-2}
(\overline{R_n(\{1\})};\Q))=-p_1.
\ees
This reflects the fact that for $n\geq 2$,
$R_n(\{1\})\setminus\{\hat{0},\hat{1}\}\cong\Pi_n\setminus\{\hat{0}\}$
is contractible. A more interesting consequence is:
\bcr \label{fibrecor}
For $n\geq 1$, $\tH_{n-2}(\overline{Q_n(G)};\Q)$ 
is isomorphic to
\[ \bigoplus_{\substack{m\geq 1\\n_1,\cdots,n_m\geq 1\\
n_1+\cdots+n_m=n}}
\Ind^{G\wr S_n}_{(G\wr S_{n_1})\times\cdots\times(G\wr S_{n_m})}
(\tH_{n_1-2}(\overline{R_{n_1}(G)};\Q)\boxtimes\cdots\boxtimes
\tH_{n_m-2}(\overline{R_{n_m}(G)};\Q)) \]
as a representation of $G\wr S_n$.
\ecr
\bpf
From Theorems \ref{hanlonthm} and \ref{secondthm} we deduce
that
\bes
\begin{split} 
1+\sum_{n\geq 1}(-1)^n 
&\ch_{G\wr S_n}(\tH_{n-2}(\overline{Q_n(G)};\Q))\\
&=(1-\sum_{n\geq 1}(-1)^n \ch_{G\wr S_n}(\tH_{n-2}
(\overline{R_n(G)};\Q)))^{-1}\\
&=\sum_{m\geq 0}(\sum_{n\geq 1}(-1)^n \ch_{G\wr S_n}(\tH_{n-2}
(\overline{R_n(G)};\Q)))^m\\
&=\sum_{\substack{m\geq 0\\n_1,\cdots,n_m\geq 1}}
(-1)^{n_1+\cdots+n_m}
\prod_{i=1}^m \ch_{G\wr S_{n_i}}(\tH_{n_i-2}
(\overline{R_{n_i}(G)};\Q)).
\end{split}
\ees
Since multiplication of Frobenius characteristics 
corresponds to induction product of
representations (\cite[Chapter I, Appendix B, (6.3)]{macdonald}), 
this gives the result.
\epf
Corollary \ref{fibrecor} can also be proved by applying 
a poset fibre theorem of Bj\"orner, Wachs, and Welker 
to the `forgetful' poset map
$Q_n(G)\to\Pi(\{0,1,\cdots,n\})$; see 
\cite[(5.3.5)]{wachs}.

For $Q_n^\sim(G)$, assuming that $G\neq\{1\}$,
we will prove the following result in Section 4:
\bth \label{thirdthm}
In $\Aa_G$ we have the equation
\[ 1+\sum_{n\geq 2}(-1)^n \ch_{G\wr S_n}(\tH_{n-2}
(\overline{Q_n^\sim(G)};\Q))=(1+\sum_{c\in G_*}\frac{|c|}{|G|}p_1(c))
(\Exp_G\circ L)^{-1}. \]
\eth
\noindent
The non-equivariant version is
\beq
1+\sum_{n\geq 2}(-1)^n \dim \tH_{n-2}
(\overline{Q_n^\sim(G)};\Q)\frac{x^n}{|G|^n n!}=
(1+\frac{x}{|G|})(1+x)^{-1/|G|},
\eeq
equivalent to the result which is well known at least for cyclic $G$
(see \cite[Corollary 6.86]{orlikterao}):
\beq
\dim \tH_{n-2}
(\overline{Q_n^\sim(G)};\Q)=(n-1)(|G|-1)(|G|+1)(2|G|+1)\cdots((n-2)|G|+1).
\eeq
A further consequence of Theorem \ref{thirdthm} is:
\bcr
For $n\geq 2$, $\tH_{n-2}(\overline{Q_n(G)};\Q)$ 
is isomorphic to
\[ \Ind^{G\wr S_n}_{G^n}(\bbone)\oplus\bigoplus_{k=0}^{n-2}
\Ind^{G\wr S_n}_{G^k\times(G\wr S_{n-k})}
(\bbone\boxtimes
\tH_{n-k-2}(\overline{Q_{n-k}^\sim(G)};\Q)) \]
as a representation of $G\wr S_n$.
\ecr
\bpf
From Theorems \ref{hanlonthm} and \ref{thirdthm} we deduce that
\bes
\begin{split}
1+&\sum_{n\geq 1}(-1)^n \ch_{G\wr S_n}(\tH_{n-2}
(\overline{Q_n(G)};\Q))\\
&=(1+\sum_{c\in G_*}\frac{|c|}{|G|}p_1(c))^{-1}
(1+\sum_{m\geq 2}(-1)^m \ch_{G\wr S_m}(\tH_{m-2}
(\overline{Q_m^\sim(G)};\Q))\\
&=1+\sum_{n\geq 1}(-1)^n\ch_{G\wr S_1}(\bbone)^n\\
&\qquad+\sum_{\substack{n\geq 2\\0\leq k\leq n-2}}
(-1)^n\ch_{G\wr S_1}(\bbone)^k\ch_{G\wr S_{n-k}}(\tH_{n-k-2}
(\overline{Q_{n-k}^\sim(G)};\Q)),
\end{split}
\ees
which implies the claim.
\epf
Perhaps this Corollary too follows from a suitable poset fibre theorem.

To state the results for the Calderbank-Hanlon-Robinson-style sub-posets we need
a bit more notation. For any $d\geq 2$, $\Exp_G^\zeromodd$ denotes 
the sum of all terms of $\Exp_G$ whose degree is $\equiv 0$ mod $d$, 
and $\Exp_G^\nonzeromodd$
denotes the sum of the other terms; we use similar notations with $0$ mod $d$
replaced by $1$ mod $d$, and with subscripts omitted when $G=\{1\}$.
Since $\Exp^\onemodd$ is an element of $\Aa_{\{1\},+}$ whose degree-$1$
term is $p_1$, it has a unique two-sided plethystic inverse in
$\Aa_{\{1\},+}$, which we write as $(\Exp^\onemodd)^{[-1]}$.

In Section 4 we will prove the following.
\bth \label{1moddthm}
In $\Aa_G$ we have the equation
\bes
\begin{split} 
1+\sum_{n\geq 1}(-1)^{\lceil\frac{n}{d}\rceil} 
&\ch_{G\wr S_n}(\tH_{\lceil\frac{n}{d}\rceil-2}
(\overline{Q_n^\onemodd(G)};\Q))\\
&=[(1-\Exp_G^\nonzeromodd)(\Exp_G^\zeromodd)^{-1}]\circ (\Exp^\onemodd)^{[-1]}.
\end{split}
\ees
\eth
\noindent
In the $d=2$ case, the right-hand side could be written more suggestively
as $(\Sech_G-\Tanh_G)\circ\Arcsinh$.
The non-equivariant version of this special case is:
\beq
\begin{split}
1+\sum_{n\geq 1}(-1)^{\lceil\frac{n}{2}\rceil} 
&\dim\tH_{\lceil\frac{n}{2}\rceil-2}
(\overline{Q_n^{1\,\mathrm{mod}\,2}(G)};\Q)
\frac{x^n}{|G|^n n!}\\
&=\sech(\frac{1}{|G|}\arcsinh(x))-\tanh(\frac{1}{|G|}\arcsinh(x)).
\end{split}
\eeq
The $G=\{\pm 1\}$ special case of this is equivalent to:
\beq
\dim\tH_{\lceil\frac{n}{2}\rceil-2}
(\overline{Q_n^{1\,\mathrm{mod}\,2}(\{\pm 1\})};\Q)
=\left\{\begin{array}{cl}
\frac{(2n)!}{2^n(n+1)!},&\text{ when $n$ is even,}\\
\frac{n!(n-1)!}{(\frac{n+1}{2})!(\frac{n-1}{2})!},
&\text{ when $n$ is odd.}
\end{array}
\right.
\eeq
Note that the $G=\{1\}$ special case of Theorem \ref{1moddthm} is known:
\bes
\begin{split} 
1+\sum_{n\geq 1}(-1)^{\lceil\frac{n}{d}\rceil} 
&\ch_{S_n}(\tH_{\lceil\frac{n}{d}\rceil-2}
(\overline{Q_n^\onemodd(\{1\})};\Q))\\
&=[(1-\Exp^\nonzeromodd)(\Exp^\zeromodd)^{-1}]\circ (\Exp^\onemodd)^{[-1]}
\end{split}
\ees
can also be obtained by applying $\frac{\partial}{\partial p_1}$
to both sides of \cite[Theorem 4.7]{chr}, which in our notation is
\beq \label{chr1modeqn}
\begin{split} 
p_1+\sum_{n\geq 2}(-1)^{\lceil\frac{n-1}{d}\rceil} 
&\ch_{S_n}(\tH_{\lceil\frac{n-1}{d}\rceil-2}
(\overline{\Pi_{n}^{(1,d)}};\Q))\\
&=(1+p_1-\Exp^{\nononemodd})\circ (\Exp^\onemodd)^{[-1]}.
\end{split}
\eeq

Finally, we have the result for the `$d$-divisible' sub-poset of the Dowling
lattice, also to be proved in Section 4.
\bth \label{0moddthm}
In $\Aa_G$ we have the equation
\bes
\begin{split}
\sum_{n\geq 1}(-1)^{\lfloor\frac{n}{d}\rfloor+1} 
&\ch_{G\wr S_n}(\tH_{\lfloor\frac{n}{d}\rfloor-1}
(\overline{Q_n^\zeromodd(G)};\Q))\\
&=1-\Exp_G\cdot(\Exp_G\circ L\circ(\Exp^\zeromodd-1))^{-1}.
\end{split}
\ees
\eth
\noindent
The non-equivariant version of the $d=2$ case is:
\beq
\begin{split}
\sum_{n\geq 1}(-1)^{\lfloor\frac{n}{2}\rfloor+1} 
\dim\tH_{\lfloor\frac{n}{2}\rfloor-1}
&(\overline{Q_n^{0\,\mathrm{mod}\,2}(G)};\Q)
\frac{x^n}{|G|^n n!}\\
&=1-(1+\tanh(x))^{1/|G|}.
\end{split}
\eeq
Note that in the $G=\{1\}$ case of Theorem \ref{0moddthm}, 
the expression $\Exp\circ L$ collapses to $1+p_1$, so we have
\bes
\sum_{n\geq 1}(-1)^{\lfloor\frac{n}{d}\rfloor+1} 
\ch_{S_n}(\tH_{\lfloor\frac{n}{d}\rfloor-1}
(\overline{Q_n^\zeromodd(\{1\})};\Q))
=1-\frac{\Exp}{\Exp^\zeromodd}.
\ees
Taking only the terms of degree $\equiv -1$ mod $d$ in this formula, we get
\bes
\sum_{\substack{n\geq d-1\\n\equiv -1\,\mathrm{mod}\,d}}(-1)^{\frac{n+1}{d}} 
\ch_{S_n}(\tH_{\frac{n+1}{d}-2}
(\overline{Q_n^\zeromodd(\{1\})};\Q))
=-\frac{\Exp^{-1\,\mathrm{mod}\,d}}{\Exp^\zeromodd}.
\ees
This can also be obtained by applying $\frac{\partial}{\partial p_1}$
to both sides of \cite[Corollary 4.7]{chr}, which in our notation is
\beq \label{chr0modeqn}
\sum_{\substack{n\geq d\\n\equiv 0\,\mathrm{mod}\,d}}(-1)^{\frac{n}{d}} 
\ch_{S_n}(\tH_{\frac{n}{d}-2}
(\overline{\Pi_{n}^{(0,d)}};\Q))
=-L\circ(\Exp^{\zeromodd}-1).
\eeq
The fact that Theorem \ref{0moddthm} for nontrivial $G$
does not simplify to the same extent perhaps implies that no Dowling
lattice analogue of \cite[Theorem 6.5]{chr} can be found.
\section{$(G\wr\Ss)$-modules}
The technical tool we will use to prove the Theorems stated in the previous
section is an extension of 
Joyal's theory of tensor species (for which see \cite{joyalsln}
or the textbook \cite{bll}) to the case
of wreath products $G\wr S_n$. This was introduced in 
\cite[Section 5]{mywreath}
(as observed at the end of that section, the assumption that 
$G$ is cyclic is unnecessary). We recall the main points here.
All vector spaces and representations are over $\Q$.

Define the category $\B_G$ whose objects
are finite sets equipped with a \textbf{free} action of $G$, and whose 
morphisms are $G$-equivariant \textbf{bijections}. Thus every object
in $\B_G$ is isomorphic to $G\times [n]$ for a unique $n\in\N$. 
The automorphism group of
$G\times [n]$ in $\B_G$ can obviously be identified with the wreath product
$G\wr S_n$.
A \textbf{tensor species}, 
also known as an \textbf{$\Ss$-module}, is a functor from
$\B_{\{1\}}$ to vector spaces;
the natural generalization is as follows.
\bdf
A \textbf{$(G\wr\Ss)$-module} (over $\Q$) is a functor $U$ from $\B_G$ to
the category of finite-dimensional vector
spaces over $\Q$. That is, to each finite set $I$ with a free $G$-action
it associates a finite-dimensional vector space $U(I)$, and to 
each $G$-equivariant bijection
$f:I\isomto J$ between such sets it associates an isomorphism
$U(f):U(I)\isomto U(J)$.
\edf
(In \cite{mywreath}, where $G$ was cyclic of order $r$, I called this
a $\B_r$-module.)

Clearly any $(G\wr\Ss)$-module $U$ gives rise to a sequence
$(U(G\times [n]))_{n\geq 0}$ of representations of the various
wreath products $G\wr S_n$; moreover, $U$ is determined up to
isomorphism (in the usual sense of isomorphism for functors) by this
sequence of representations. This means that $U$ is determined up to 
isomorphism by its \textbf{character}
\beq
\ch(U):= \sum_{n\geq 0}\ch_{G\wr S_n}(U(G\times [n]))\in\Aa_G.
\eeq
The convenience of 
defining $U$ as a functor rather
than just a sequence of representations will become clear shortly.

An important example is the trivial $(G\wr\Ss)$-module $\bbone_G$, which
is defined by $\bbone_G(I)=\Q$ for all objects $I$ of $\B_G$, and
$\bbone_G(f)=\mathrm{id}$ for all morphisms $f$ of $\B_G$. Clearly
$\bbone_G(G\times [n])$ is the trivial representation of $G\wr S_n$, so
$\ch(\bbone_G)=\Exp_G$.

For any $(G\wr\Ss)$-module $U$ we can define various sub-$(G\wr\Ss)$-modules
by imposing a restriction on degree, which we will write as a superscript.
For instance, $U^{\onemodd}$ is the $(G\wr\Ss)$-module defined by
$U^{\onemodd}(I)=U(I)$ when $\frac{|I|}{|G|}\equiv 1$ mod $d$, 
and $U^{\onemodd}(I)=0$
when $\frac{|I|}{|G|}\not\equiv 1$ mod $d$; the definition of
$U^{\onemodd}$ on morphisms of $\B_G$ is the same as that of $U$ when
the morphisms are of the form $f:I\isomto J$ for $\frac{|I|}{|G|}
=\frac{|J|}{|G|}\equiv 1$ mod $d$,
and zero otherwise. Clearly $\ch(U^{\onemodd})$ is the sum of all terms of
$\ch(U)$ whose degree is $\equiv 1$ mod $d$. Similarly we define
$U^{\zeromodd}$, $U^{\nononemodd}$, $U^{\nonzeromodd}$, and $U^{\geq 1}$.

Now the main point of Joyal's theory is that certain natural operations
of tensor species, called sum, product, and substitution
(or partitional composition), correspond to the analogous operations 
on their characters;
the analogue of substitution is plethysm. These
operations can be extended to our context as follows.
We define the sum and product of two $(G\wr\Ss)$-modules $U$ and $V$ by
\begin{equation} \label{sumproddefeqn}
\begin{split}
(U+V)(I)&=U(I)\oplus V(I),\\
(U\cdot V)(I)&=\bigoplus_{\substack{J\subseteq I\\G\text{-stable}}}
U(J)\otimes V(I\setminus J),
\end{split}
\end{equation}
for any object $I$ of $\B_G$; the definition on morphisms is the obvious one.
(It is clear that if $J$ is a $G$-stable subset of an object of $\B_G$, then
both $J$ and $I\setminus J$ are objects of $\B_G$.)

If $U$ is an $\Ss$-module and $V$ is a $(G\wr\Ss)$-module such that
$V(\emptyset)=0$, we define the
substitution $U\circ V$, a $(G\wr\Ss)$-module, by
\begin{equation}
(U\circ V)(I)=\bigoplus_{\substack{\pi\in\Pi(I)\\\text{partwise $G$-stable}}}
\left(U(\pi)\otimes\bigotimes_{J\in\pi}V(J)\right)
\end{equation}
for any object $I$ of $\B_G$.
Here $\Pi(I)$ is the set of set partitions of $I$ (a set partition of $I$
is a set of nonempty disjoint subsets whose union is $I$), and $\pi\in\Pi(I)$
is `partwise $G$-stable' if $g.J=J$ for all $g\in G$ and all 
parts $J\in\pi$. The definition
on morphisms is the obvious one.

If $U$ is a $(G\wr\Ss)$-module and $V$ is an $\Ss$-module
such that $V(\emptyset)=0$, we define the
substitution $U\circ V$, a $(G\wr\Ss)$-module, by
\begin{equation} 
(U\circ V)(I)=\bigoplus_{\substack{\pi\in\Pi(I)\\\pi\in\B_G}}
\left(U(\pi)\otimes\bigotimes_{\CalO\in G\setminus\pi}V(\CalO)\right)
\end{equation}
for any object $I$ of $\B_G$.
Here the condition $\pi\in\B_G$ means just that $G$ acts freely on
the set of parts of the partition, i.e.\ for all $1\neq g\in G$ and
$J\in\pi$, $J\neq g.J\in\pi$. If $\CalO$ is a $G$-orbit on the
set of parts, $V(\CalO)$ should be thought of as $V(J)$ for some
$J\in\CalO$, the choice making no difference up to isomorphism; but in
order to be able to repeat the mantra that ``the definition on morphisms
is the obvious one'', we must make the more canonical definition that
\[ V(\CalO):=\{(v_J)\in\prod_{J\in\CalO} V(J)\,|\,
v_{g.J}=V(g|_J)(v_J),\ \forall J\in\CalO, g\in G\}. \]
By functoriality and freeness, the choice of one $v_J$ uniquely determines
the whole $\CalO$-tuple.

We have the following generalization of Joyal's result (which is
the case $G=\{1\}$).
\bth
\ben
\item If $U$ and $V$ are $(G\wr\Ss)$-modules,
\[ \ch(U+ V)=\ch(U)+\ch(V),\ \ch(U\cdot V)=\ch(U)\ch(V). \]
\item If $U$ is an $\Ss$-module and $V$ a $(G\wr\Ss)$-module such that
$V(\emptyset)=0$,
\[ \ch(U\circ V)=\ch(U)\circ\ch(V). \]
\item If $U$ is a $(G\wr\Ss)$-module and $V$ an $\Ss$-module such that
$V(\emptyset)=0$,
\[ \ch(U\circ V)=\ch(U)\circ\ch(V). \]
\een
\eth
\bpf
Probably any proof of Joyal's result could be modified to prove this;
see \cite[Proposition 5.1, Theorems 5.6 and 5.9]{mywreath} 
for such a modification
of the `analytic functors' proof (the assumption there that $G$ was cyclic
was never actually used, except in the trivial sense that conjugacy classes
were replaced by elements in the notation). The key idea is to associate
to each $(G\wr\Ss)$-module $U$ a functor
\beq
F_U:M\mapsto\bigoplus_{n\geq 0}(U(G\times [n])\otimes M^{\otimes n})^{G\wr S_n}
\eeq
from representations of $G$ to vector spaces. One then shows that the left
plethystic action corresponds to \textbf{post}composing such a functor with an
analytic functor from vector spaces to vector spaces, whereas the right
plethystic action corresponds to \textbf{pre}composing such a functor with
the functor from representations of $G$ to representations of $G$
\textbf{induced} by an analytic functor from vector spaces to vector spaces.
\epf

Since the representations we want to apply this to are homology
groups of posets, we will need to have a `super' version of the above.
So we define a \textbf{super-$(G\wr\Ss)$-module} to be a functor from $\B_G$
to the category of finite-dimensional $\Z/2\Z$-graded vector spaces over
$\Q$. If $U$ is a super-$(G\wr\Ss)$-module, define its super-character
\beq
\sch(U):= \sum_{n\geq 0}\ch_{G\wr S_n}(U(G\times [n])_{\bar{0}})-
\ch_{G\wr S_n}(U(G\times [n])_{\bar{1}}).
\eeq
Any $(G\wr\Ss)$-module $U$ may be viewed as a super-$(G\wr\Ss)$-module
which is purely even, so that $\sch(U)=\ch(U)$.
The above definitions of sum, product, and substitution
can be carried over to the super context,
using the usual sign-commutativity convention for tensor products.
\bth \label{superthm}
\ben
\item If $U$ and $V$ are super-$(G\wr\Ss)$-modules,
\[ \sch(U+ V)=\sch(U)+\sch(V),\ \sch(U\cdot V)=\sch(U)\sch(V). \]
\item If $U$ is a super-$\Ss$-module and $V$ a 
super-$(G\wr\Ss)$-module such that
$V(\emptyset)=0$,
\[ \sch(U\circ V)=\sch(U)\circ\sch(V). \]
\item If $U$ is a super-$(G\wr\Ss)$-module and $V$ a 
super-$\Ss$-module such that
$V(\emptyset)=0$,
\[ \sch(U\circ V)=\sch(U)\circ\sch(V). \]
\een
\eth
\bpf
A more complicated version of this theorem was proved in 
\cite[Section 7]{mywreath}. To deduce the present version, set $q\to 1$
in \cite[Corollaries 7.3 and 7.6]{mywreath} to obtain (1) and (3);
the analogue of (2) was not stated there but follows by the same method.
\epf
\section{Proof of the results}
In this section we prove Theorems \ref{hanlonthm}, \ref{secondthm},
\ref{thirdthm},
\ref{1moddthm} and \ref{0moddthm}.
To make our arguments more legible, we need a notational convention:
for any Cohen-Macaulay poset $P$ and elements $x<y$, we write
$\tH_P(x,y)$ for the top-degree reduced homology $\tH_{\ell(x,y)-2}((x,y);\Q)$
of the open interval $(x,y)\subset P$.
We view this as a super vector space of parity equal to that of $\ell(x,y)$.
Since the top degree is the only one in which the reduced homology 
could be nonzero, we have 
\beq 
\sdim \tH_P(x,y)=(-1)^{\ell(x,y)}\dim \tH_P(x,y)
=\widetilde{\chi}((x,y))=\mu_P(x,y).
\eeq
We also define $\tH_P(x,x)$ to be a one-dimensional even super vector space,
so that $\sdim\tH_P(x,x)=1=\mu_P(x,x)$.
The key fact we use, an application of the Euler-Poincar\'e principle
observed by Sundaram (see \cite[Section 1]{sundaram} and
\cite[Theorem 4.4.1]{wachs}), is as follows:
\bth \label{sundaramthm}
If $P$ is a Cohen-Macaulay pure bounded poset where 
$\hat{0}\neq\hat{1}$ \textup{(}i.e.\ $P$ has more than one element\textup{)},
then $\bigoplus_{x\in P} \tH_P(\hat{0},x)$ is balanced; that is, its
odd and even parts are isomorphic, as vector spaces and as representations
of any group that acts on $P$. The same is true for
$\bigoplus_{x\in P} \tH_P(x,\hat{1})$.
\eth
\noindent
Of course the second statement is just the first applied to the dual poset.
(Taking sdim, we recover the familiar recursive properties
of the M\"obius function.) 

In order to apply the theory of the previous section, we need to define
our posets functorially. We illustrate by rewriting Sundaram's proof
of Stanley's theorem \eqref{stanleyeqn} on 
the partition lattices (see \cite[Theorem 4.4.7]{wachs}).
We have already defined $\Pi(I)$ as the lattice
of set partitions of the set $I$; with an obvious definition on morphisms,
this constitutes a functor $\Pi$ from $\B_{\{1\}}$ to the category of
posets. We then define three related 
super-$\Ss$-modules, $\tH_\Pi$, $\WH_\Pi$, and $\WH_\Pi^*$. 
The definitions on objects
of $\B_{\{1\}}$ are:
\bes
\begin{split} 
\tH_\Pi(I)&=\tH_{\Pi(I)}(\hat{0},\hat{1}),\\
\WH_\Pi(I)&=\bigoplus_{\pi\in\Pi(I)}\tH_{\Pi(I)}(\hat{0},\pi),
\text{ and }\\
\WH_\Pi^*(I)&=\bigoplus_{\pi\in\Pi(I)}\tH_{\Pi(I)}(\pi,\hat{1}),
\end{split}
\ees
and the definitions on morphisms are the obvious one.
(By our convention that $\Pi(\emptyset)$ has one element,
$\tH_\Pi(\emptyset)$, $\WH_\Pi(\emptyset)$, and $\WH_\Pi^*(\emptyset)$ 
are one-dimensional and of even parity.)
Sundaram's argument rests on the recursive property of partition lattices,
namely that for any $\pi\in\Pi_n$, the principal upper order ideal
$[\pi,\hat{1}]$ is isomorphic to $\Pi_{|\pi|}$. With functorial
language we can state this more precisely: for any $\pi\in\Pi(I)$, 
$[\pi,\hat{1}]$ is \textbf{canonically} isomorphic to $\Pi(\pi)$.
Recalling the definition of substitution from the previous section,
we see that we have an isomorphism of super-$\Ss$-modules:
\beq \label{firstpartitioneqn}
\WH_\Pi^*\cong \tH_\Pi\circ\bbone_{\{1\}}^{\geq 1}.
\eeq
Taking $\sch$ and applying
the $G=\{1\}$ case of Theorem \ref{superthm}, we get
\bes
\sch(\WH_\Pi^*)=\sch(\tH_\Pi)\circ(\Exp-1).
\ees
Now by Theorem \ref{sundaramthm}, $\sch(\WH_\Pi^*)=1+p_1$ (only the
one-element posets contribute).
We deduce that $\sch(\tH_\Pi)=1+L$, which is exactly \eqref{stanleyeqn}.
A slightly more complicated argument uses lower order ideals: for any
$\pi\in\Pi(I)$, 
$[\hat{0},\pi]$ is canonically isomorphic to
$\prod_{K\in\pi}\Pi(K)$. Applying the K\"unneth formula
\cite[second statement of Theorem 5.1.5]{wachs}, 
we see that we have an isomorphism of super-$\Ss$-modules:
\beq \label{secondpartitioneqn}
\WH_\Pi\cong \bbone_{\{1\}}\circ\tH_\Pi^{\geq 1},
\eeq
and the result follows as before. Notice how the sign convention in
the definition of substitution of super-$\Ss$-modules takes into
account the sign-commutativity of the K\"unneth formula.

With these arguments as models, we turn to our sub-posets of the Dowling
lattices. First we have to reinterpret them as functors 
$Q$, $R$, $Q^\sim$, $Q^\onemodd$, $Q^\zeromodd$ from $\B_G$ to
the category of posets, so that $Q(G\times [n])=Q_n(G)$,
$R(G\times [n])=R_n(G)$, and so on. The definitions in Section 1 were
deliberately written so that this is simply a matter of replacing $G\times [n]$
with a general object $I$ of $\B_G$. (Following usual conventions, this
will result in $Q(\emptyset)$, $R(\emptyset)$, $Q^\sim(\emptyset)$,
$Q^\onemodd(\emptyset)$,
and $Q^\zeromodd(\emptyset)$ all being one-element posets; we also stipulate
that when $\frac{|J|}{|G|}=1$, $Q^\sim(J)$ is a one-element poset.)
We then define three
super-$(G\wr\Ss)$-modules attached to each functor, as with $\Pi$:
$\tH_Q$, $\WH_Q$, $\WH_Q^*$, $\tH_R$, $\WH_R$, and $\WH_R^*$, and so
forth.

\noindent
\bpf \textbf{(Theorem \ref{hanlonthm})}
For any $(J,\pi)\in Q(I)$,
$[(J,\pi),\hat{1}]$ is isomorphic to $Q(\pi)$, so
\[ \WH_Q^*(I)\cong\bigoplus_{\substack{J\subseteq I\\\text{$G$-stable}}}
\bigoplus_{\substack{\pi\in\Pi(I\setminus J)\\\pi\in\B_G}}
\tH_Q(\pi). \]
Clearly this amounts to
an isomorphism of super-$(G\wr\Ss)$-modules:
\beq
\WH_Q^*\cong \bbone_G\cdot(\tH_Q\circ\bbone_{\{1\}}^{\geq 1}).
\eeq
Taking $\sch$ and applying parts (1) and (3) of Theorem \ref{superthm},
we get
\bes
\sch(\WH_Q^*)=\Exp_G\cdot(\sch(\tH_Q)\circ(\Exp-1)).
\ees
Now by Theorem \ref{sundaramthm}, $\sch(\WH_Q^*)=1$ (only the one-element
poset $Q(\emptyset)$ contributes). Hence
\beq
\sch(\tH_Q)=\Exp_G^{-1}\circ L=(\Exp_G\circ L)^{-1},
\eeq
which is exactly the statement. For reference, we also give the alternative
proof using $\WH_Q$. For any $(J,\pi)\in Q(I)$,
$[\hat{0},(J,\pi)]$ is isomorphic to 
$Q(J)\times\prod_{\CalO\in G\setminus\pi}\Pi(K_{\CalO})$ where
$K_{\CalO}$ denotes a representative of the orbit $\CalO$,
so by the same K\"unneth
formula as above,
\[ \WH_Q(I)\cong\bigoplus_{\substack{J\subseteq I\\\text{$G$-stable}}}
\tH_Q(J)\otimes
(\bigoplus_{\substack{\pi\in\Pi(I\setminus J)\\\pi\in\B_G}}
\bigotimes_{\CalO\in G\setminus\pi} \tH_\Pi(K_{\CalO})). \]
Clearly this amounts to
an isomorphism of super-$(G\wr\Ss)$-modules:
\beq \label{whitneyeqn}
\WH_Q\cong \tH_Q\cdot(\bbone_G\circ \tH_\Pi^{\geq 1}).
\eeq
Using $\sch(\WH_Q)=1$ and $\sch(\tH_\Pi^{\geq 1})=L$, we
reach the result again. 
\epf

\noindent
\bpf \textbf{(Theorem \ref{secondthm})}
For any
$(\emptyset,\pi)\in R(I)$, $[(\emptyset,\pi),\hat{1}]$
is isomorphic to $R(\pi)$. Hence we have
an isomorphism of super-$(G\wr\Ss)$-modules:
\beq
\WH_R^*\cong \bbone_G^{\geq 1}+\tH_R\circ\bbone_{\{1\}}^{\geq 1},
\eeq
where the first term comes from the $\tH_{R(I)}(\hat{1},\hat{1})$ terms.
Taking $\sch$ and applying parts (1) and (3) of Theorem \ref{superthm},
we get
\bes
\sch(\WH_R^*)=\Exp_G-1+\sch(\tH_R)\circ(\Exp-1).
\ees
Now by Theorem \ref{sundaramthm}, $\sch(\WH_R^*)=1$ (only the one-element
poset $R(\emptyset)$ contributes). Hence
\beq
\sch(\tH_R)=(2-\Exp_G)\circ L=2-\Exp_G\circ L,
\eeq
which (subtracting $1$ from both sides) is exactly the statement.
The alternative proof would use the isomorphism
\beq \label{secondwhitneyeqn}
\WH_R\cong\tH_R^{\geq 1}+\bbone_G\circ\tH_\Pi^{\geq 1}
\eeq
of super-$(G\wr\Ss)$-modules.
\epf

\noindent
\bpf \textbf{(Theorem \ref{thirdthm})}
Here the proof via $\WH_{Q^\sim}$ is more convenient. 
For any $(J,\pi)\in Q^\sim(I)$,
$[\hat{0},(J,\pi)]\cong
Q^\sim(J)\times\prod_{\CalO\in G\setminus\pi}\Pi(K_{\CalO})$, so
the analogue of
\eqref{whitneyeqn} is
\beq \label{thirdwhitneyeqn}
\WH_{Q^\sim}\cong \tH_{Q^\sim}^{\neq 1}
\cdot(\bbone_G\circ \tH_\Pi^{\geq 1}).
\eeq
Taking $\sch$ and applying parts (1) and (3) of Theorem \ref{superthm},
we get
\bes
\sch(\WH_{Q^\sim})=\sch(\tH_{Q^\sim})^{\neq 1}\cdot(\Exp_G\circ L).
\ees
Now by Theorem \ref{sundaramthm}, $\sch(\WH_{Q^\sim})=
1+\sum_{c\in G_*}\frac{|c|}{|G|}p_1(c)$ (only the one-element
posets $Q^\sim(\emptyset)$ and $Q^\sim(G\times [1])$ contribute). Hence
\beq
\sch(\tH_{Q^\sim})^{\neq 1}=(1+\sum_{c\in G_*}\frac{|c|}{|G|}p_1(c))
(\Exp_G\circ L)^{-1},
\eeq
which is exactly the statement.
\epf

\noindent
\bpf \textbf{(Theorem \ref{1moddthm})}
For any
$(J,\pi)\in Q^\onemodd(I)$, $[(J,\pi),\hat{1}]$
is isomorphic to $Q^\onemodd(\pi)$. So we have
an isomorphism of super-$(G\wr\Ss)$-modules:
\beq
\WH_{Q^\onemodd}^*\cong \bbone_G^{\nonzeromodd}+
\bbone_G^{\zeromodd}\cdot(\tH_{Q^{\onemodd}}\circ\bbone^{\onemodd}),
\eeq
where the two terms cover the cases of $\frac{|J|}{|G|}\not\equiv 0$ mod $d$
(forcing $J=I$)
and $\frac{|J|}{|G|}\equiv 0$ mod $d$ respectively.
Taking $\sch$ as usual,
we get
\bes
\sch(\WH_{Q^\onemodd}^*)=\Exp_G^\nonzeromodd+\Exp_G^\zeromodd
(\sch(\tH_{Q^\onemodd})\circ\Exp^\onemodd).
\ees
Now by Theorem \ref{sundaramthm}, $\sch(\WH_{Q^\onemodd}^*)=1$. Hence
\beq
\sch(\tH_{Q^\onemodd})=
[(1-\Exp_G^\nonzeromodd)(\Exp_G^\zeromodd)^{-1}]\circ (\Exp^\onemodd)^{[-1]},
\eeq
which is exactly the statement. The alternative proof would use the
isomorphism
\beq \label{fourthwhitneyeqn}
\WH_{Q^\onemodd}\cong\tH_{Q^\onemodd}^{\nonzeromodd}
+\tH_{Q^\onemodd}^{\zeromodd}\cdot
(\bbone_G\circ\tH_{\Pi^{\onemodd}}^{\onemodd}),
\eeq
and the Calderbank-Hanlon-Robinson result \eqref{chr1modeqn},
which implies that $\sch(\tH_{\Pi^{\onemodd}}^{\onemodd})=
(\Exp^\onemodd)^{[-1]}$. (Here $\Pi^{\onemodd}$ denotes the functor
from $\B_{\{1\}}$ to the category of posets such that
$\Pi^{\onemodd}([n])$ is the poset $\Pi_n^{(1,d)}$ studied in \cite{chr}.)
\epf

\noindent
\bpf \textbf{(Theorem \ref{0moddthm})}
For any
$(J,\pi)\in Q^\zeromodd(I)$ except $\hat{0}$, $[(J,\pi),\hat{1}]$
is isomorphic to $Q(\pi)$. So we have
an isomorphism of super-$(G\wr\Ss)$-modules:
\beq
\WH_{Q^\zeromodd}^*\cong\tH_{Q^\zeromodd}^{\geq 1}+
\bbone_G\cdot(\tH_{Q}\circ\bbone^{\zeromodd,\geq d}),
\eeq
where the two terms cover respectively
the cases where $(J,\pi)=\hat{0}$,
$I\neq\emptyset$,
and where $|K|\equiv 0$ mod $d$ for all $K\in\pi$.
Taking $\sch$ and using Theorem \ref{hanlonthm},
we get
\bes
\begin{split}
\sch(\WH_{Q^\zeromodd}^*)&=\sch(\tH_{Q^\zeromodd})-1\\
&\qquad+\Exp_G\cdot((\Exp_G\circ L)^{-1}\circ(\Exp^\zeromodd-1)).
\end{split}
\ees
Now by Theorem \ref{sundaramthm}, $\sch(\WH_{Q^\zeromodd}^*)=1$. Hence
\beq \label{zeromodresulteqn}
\sch(\tH_{Q^\zeromodd})=
2-\Exp_G\cdot(\Exp_G\circ L\circ(\Exp^\zeromodd-1))^{-1},
\eeq
which (subtracting $1$ from both sides) is exactly the statement.
The alternative proof via $\WH_{Q^\zeromodd}$ requires a bit of care,
since for $\hat{0}\neq(J,\pi)\in Q^\zeromodd(I)$, it is not the closed
interval $[\hat{0},(J,\pi)]$ but rather the semi-closed interval
$(\hat{0},(J,\pi)]$ which is naturally a product of smaller posets.
Hence one must use
the `once-suspended' K\"unneth formula
\cite[first statement of Theorem 5.1.5]{wachs}. The upshot is
the following isomorphism of super-$(G\wr\Ss)$-modules, where
$U[1]$ denotes $U$ with parities interchanged:
\beq \label{suspendedeqn}
\begin{split}
\WH_{Q^\zeromodd}[1]&\cong\bbone_G[1]+\bbone_G^{\geq 1}\circ
(\tH_{\Pi^\zeromodd}^{\zeromodd,\geq d}[1])\\
&\quad+\tH_{Q^\zeromodd}^{\geq 1}[1]
\cdot(\bbone_G\circ\tH_{\Pi^\zeromodd}^{\zeromodd,\geq d}[1]).
\end{split}
\eeq
Here the first term corresponds to the minimum element $\hat{0}$,
the second to the non-minimum elements of the form $(\emptyset,\pi)$, 
and the third to the elements of the form
$(J,\pi)$ with $J\neq\emptyset$; $\Pi^{\zeromodd}$ denotes the functor
from $\B_{\{1\}}$ to the category of posets such that
$\Pi^{\zeromodd}([n])$ is the poset $\Pi_n^{(0,d)}$ studied in \cite{chr}.
Now the Calderbank-Hanlon-Robinson result \eqref{chr0modeqn} implies
$\sch(\tH_{\Pi^\zeromodd}^{\zeromodd,\geq d}[1])=
L\circ(\Exp^\zeromodd-1)$, so we obtain
\bes
\begin{split}
-1&=-\Exp_G+(\Exp_G-1)\circ L\circ (\Exp^\zeromodd-1)\\
&\quad+(1-\sch(\tH_{Q^\zeromodd}))(\Exp_G\circ L\circ (\Exp^\zeromodd-1)),
\end{split}
\ees
which simplifies to \eqref{zeromodresulteqn} again.
\epf
\section{Whitney homology}
The name $\WH$ for the super-$(G\wr\Ss)$-modules used in the previous
section of course stands for `Whitney homology'; but properly speaking,
the Whitney homology of a poset has a $\Z$-grading, not a 
$(\Z/2\Z)$-grading. Recall that if $P$ is a Cohen-Macaulay poset with minimum
element $\hat{0}$, its
(rational) Whitney homology groups are defined by
\[ \WH_i(P):=\bigoplus_{\substack{x\in P\\\rk(x)=i}}\tH_P(\hat{0},x),
\text{ for all $i\in\Z$.} \]
If a group $\Gamma$ acts on the poset $P$, it also acts on each
$\WH_i(P)$, and the characters of these Whitney homology representations
encapsulate the `equivariant characteristic polynomials':
\[ \sum_{i\in\Z}\tr(\gamma,\WH_i(P))\,(-t)^i
=\sum_{x\in P^\gamma}\mu_{P^\gamma}(\hat{0},x)\,t^{\rk_P(x)},
\text{ for all $\gamma\in\Gamma$.} \]
It is an easy matter to find formulae analogous to those in \S2
for these characters in the case of of our posets,
since we have effectively already worked out the relationship
between the Whitney homology and its highest-degree part.

We need to introduce the concept of a \textbf{graded $(G\wr\Ss)$-module},
which is merely a functor from $\B_G$ to the category of 
finite-dimensional $\Z$-graded vector spaces over $\Q$. 
Any $(G\wr\Ss)$-module may be regarded as a graded
$(G\wr\Ss)$-module concentrated in degree $0$. If $U$ is a
graded $(G\wr\Ss)$-module, we define its graded character to be
\[ \ch_t(U):=\sum_{n\geq 0}\sum_{i\in\Z}
\ch_{G\wr S_n}(U(G\times [n])_i)\,(-t)^i, \]
an element of $\Aa_G\otimes_\Q \Q[t,t^{-1}]$. 
Note that setting $t\to 1$ recovers the super-character of $U$ regarded
as a super-$(G\wr\Ss)$-module. The definitions of sum, product, and
substitution carry over to this graded context in the usual way,
incorporating the sign-commutativity in tensor products. We can also
extend the definitions of our plethystic actions to
$\Aa_G\otimes_\Q \Q[t,t^{-1}]$ by adding the rules that
$p_i\circ t=t^i$, $p_i(c)\circ t=t^i$.
\bth \label{gradedthm}
\ben
\item If $U$ and $V$ are graded $(G\wr\Ss)$-modules,
\[ \ch_t(U+ V)=\ch_t(U)+\ch_t(V),\ \ch_t(U\cdot V)=\ch_t(U)\ch_t(V). \]
\item If $U$ is a graded $\Ss$-module and $V$ a 
graded $(G\wr\Ss)$-module such that
$V(\emptyset)=0$,
\[ \ch_t(U\circ V)=\ch_t(U)\circ\ch_t(V). \]
\item If $U$ is a graded $(G\wr\Ss)$-module and $V$ a 
graded $\Ss$-module such that
$V(\emptyset)=0$,
\[ \ch_t(U\circ V)=\ch_t(U)\circ\ch_t(V). \]
\een
\eth
\bpf
Again, a more complicated version was proved in
\cite[Section 7]{mywreath}. For (1), set $q\to 1$
in \cite[Proposition 7.2]{mywreath}; for (3), set $q\to 1$
in \cite[Theorem 7.5]{mywreath}; part (2)
follows by the same method.
\epf

To apply this, we must reinterpret $\tH$, $\WH$, and $\WH^*$
as graded $(G\wr\Ss)$-modules for each of our poset functors, by viewing
each $\tH_P(x,y)$ as a graded vector space concentrated in 
degree $\ell(x,y)$. (Recall that the actual homological degree, when
$x<y$, is $\ell(x,y)-2$; the `double suspension' is required by
the K\"unneth formula \cite[second statement of Theorem 5.15]{wachs}.)
It is easy to check that every isomorphism of super-$(G\wr\Ss)$-modules
stated in the previous section remains true verbatim as an isomorphism
of graded $(G\wr\Ss)$-modules, to which we can apply $\ch_t$ and use
the appropriate parts of Theorem \ref{gradedthm}.
(The notation $[1]$ in \eqref{suspendedeqn} now denotes a shift
of grading, such that $U[1]_i=U_{i+1}$.)

The flow of information is reversed from that in the previous
section: initially we knew $\sch(\WH)$ by Theorem \ref{sundaramthm}
and deduced $\sch(\tH)$; now we can easily obtain $\ch_t(\tH)$ from this,
and deduce $\ch_t(\WH)$. To illustrate the procedure on the partition
lattice once more, the first step is to complete the second equation
of the following analogy:
\bes
\begin{split}
\sch(\tH_\Pi)&=1+\sum_{n\geq 1}(-1)^{n-1}
\ch_{S_n}(\tH_{\Pi_n}(\hat{0},\hat{1}))=1+L,\\
\ch_t(\tH_\Pi)&=1+\sum_{n\geq 1}\ch_{S_n}(\tH_{\Pi_n}(\hat{0},\hat{1}))\,
(-t)^{n-1}=\text{ ???}
\end{split}
\ees
Since the only difference from first equation to second is that the degree-$n$ term for each $n\geq 1$ is multiplied by
$t^{n-1}$ (the exponent being $l(\Pi_n)$), the answer is clearly that
\beq \label{whitneyparteqn}
\ch_t(\tH_\Pi)=1+t^{-1}L\circ tp_1.
\eeq
Hence \eqref{firstpartitioneqn} and \eqref{secondpartitioneqn},
viewed as isomorphisms of graded $\Ss$-modules, imply that
\bes
\begin{split}
\ch_t(\WH_\Pi)&=\Exp\circ t^{-1}L\circ tp_1,\\
\ch_t(\WH_\Pi^*)&=1+t^{-1}L\circ t(\Exp-1).
\end{split}
\ees
Unravelling the definitions to rewrite the left-hand sides, 
these become the known facts
\bes
\begin{split}
1+\sum_{n\geq 1}\sum_{i\in\Z}\ch_{S_n}(\WH_i(\Pi_n))\,(-t)^i&=\Exp\circ t^{-1}L\circ tp_1,\\
1+\sum_{n\geq 1}\sum_{i\in\Z}\ch_{S_n}(\WH_i(\Pi_n^*))\,(-t)^i
&=1+t^{-1}L\circ t(\Exp-1),
\end{split}
\ees
which are equivariant version of the familiar generating
functions for the characteristic polynomials of the partition lattices
and their duals.

In the remainder of the section we give
the results of applying this same procedure to our
wreath product posets. (We omit the dual forms, which
seem less interesting.)
\bth \textup{(}Hanlon\textup{)}
In $\Aa_G\otimes_\Q \Q[t,t^{-1}]$ we have the equation
\[ 1+\sum_{n\geq 1}\sum_{i\in\Z}\ch_{G\wr S_n}(\WH_i(Q_n(G)))\,(-t)^i
=\Exp_G\circ(t^{-1}-1)L\circ tp_1. \]
\eth
\bpf
From Theorem \ref{hanlonthm} we deduce that 
$\ch_t(\tH_Q)=(\Exp_G\circ L\circ tp_1)^{-1}$. Substituting this fact and
\eqref{whitneyparteqn} into \eqref{whitneyeqn} gives
\[ \ch_t(\WH_Q)=(\Exp_G\circ L\circ tp_1)^{-1}
(\Exp_G\circ t^{-1}L\circ tp_1), \]
which is the statement.
\epf
Note that this Theorem is just a rephrasing of 
\cite[Corollary 2.3]{hanlonwreath}. Its non-equivariant version is
\beq
1+\sum_{n\geq 1}\sum_{i\in\Z}\dim\WH_i(Q_n(G))\,(-t)^i\frac{x^n}{|G|^n n!}
=(1+tx)^{\frac{1}{|G|}(t^{-1}-1)},
\eeq
which is equivalent to the well-known formula
\beq
\begin{split}
\sum_{i\in\Z}&\dim\WH_i(Q_n(G))\,(-t)^i\\
&=(1-t)(1-(|G|+1)t)\cdots(1-((n-1)|G|+1)t).
\end{split}
\eeq
A similar proof using \eqref{secondwhitneyeqn} gives:
\bth
In $\Aa_G\otimes_\Q \Q[t,t^{-1}]$ we have the equation
\[ \sum_{n\geq 1}\sum_{i\in\Z}\ch_{G\wr S_n}(\WH_i(R_n(G)))\,(-t)^i
=\Exp_G\circ t^{-1}L\circ tp_1-\Exp_G\circ L\circ tp_1. \]
\eth
\noindent
The non-equivariant version is
\beq
\sum_{n\geq 1}\sum_{i\in\Z}\dim\WH_i(R_n(G))\,(-t)^i\frac{x^n}{|G|^n n!}
=(1+tx)^{\frac{1}{|G|}t^{-1}}-(1+tx)^{\frac{1}{|G|}},
\eeq
which is equivalent to
\beq
\begin{split}
\sum_{i\in\Z}&\dim\WH_i(R_n(G))\,(-t)^i\\
&=(1-|G|t)(1-2|G|t)\cdots(1-(n-1)|G|t)\\
&\qquad\qquad+(|G|-1)(2|G|-1)\cdots((n-1)|G|-1)(-t)^n.
\end{split}
\eeq
Another similar proof using \eqref{thirdwhitneyeqn} gives
(assuming $G\neq\{1\}$):
\bth
In $\Aa_G\otimes_\Q \Q[t,t^{-1}]$ we have the equation
\bes
\begin{split}
1+\sum_{n\geq 2}
&\sum_{i\in\Z}\ch_{G\wr S_n}(\WH_i(Q_n^\sim(G)))\,(-t)^i\\
&=(1+\sum_{c\in G_*}\frac{|c|}{|G|}tp_1(c))
(\Exp_G\circ (t^{-1}-1)L\circ tp_1).
\end{split}
\ees
\eth
\noindent
The non-equivariant version is
\beq
1+\sum_{n\geq 2}
\sum_{i\in\Z}\dim\WH_i(Q_n^\sim(G))\,(-t)^i\frac{x^n}{|G|^n n!}
=(1+\frac{1}{|G|}tx)
(1+tx)^{\frac{1}{|G|}(t^{-1}-1)},
\eeq
equivalent to the formula
\beq
\begin{split}
\sum_{i\in\Z}&\dim\WH_i(Q_n^\sim(G))\,(-t)^i=(1-t)(1-(|G|+1)t)\times\cdots\\
&\qquad\cdots\times(1-((n-2)|G|+1)t)
(1-(n-1)(|G|-1)t),
\end{split}
\eeq
which for cyclic $G$ is a consequence of \cite[Corollary 6.86]{orlikterao}.

To state the analogous results for $Q^\onemodd$ and $Q^\zeromodd$,
we need to abuse notation slightly. For instance, to apply
$\ch_t$ to \eqref{fourthwhitneyeqn} we need the graded version
of the equation $\sch(\tH_{\Pi^{\onemodd}}^\onemodd)=
(\Exp^\onemodd)^{[-1]}$, namely
\beq
\ch_t(\tH_{\Pi^{\onemodd}}^\onemodd)=t^{-1/d}(\Exp^\onemodd)^{[-1]}\circ
t^{1/d}p_1.
\eeq
The right-hand side makes sense because every term of
$(\Exp^\onemodd)^{[-1]}$ has degree $\equiv 1$ mod $d$, so the exponents
of $t$ all come out as integers. (Explicitly, the degree-$n$ term is multiplied by $t^{\frac{n-1}{d}}$, as required since $l(\Pi_n^{(1,d)})=
\frac{n-1}{d}$ when $n\equiv 1$ mod $d$).
Similar remarks apply to the right-hand side of the following result,
deduced from \eqref{fourthwhitneyeqn}.
\bth \label{mainthm}
In $\Aa_G\otimes_\Q \Q[t,t^{-1}]$ we have the equation
\bes
\begin{split} 
1+\sum_{n\geq 1}
&\sum_{i\in\Z}\ch_{G\wr S_n}(\WH_i(Q_n^\onemodd(G)))\,(-t)^i\\
&=-\sum_{j=1}^{d-1}t^{\frac{d-j}{d}}(\Exp_G^{j\,\mathrm{mod}\,d}
(\Exp_G^\zeromodd)^{-1})\circ(\Exp^\onemodd)^{[-1]}\circ t^{1/d}p_1\\
&\qquad+(\Exp_G^\zeromodd\circ(\Exp^\onemodd)^{[-1]}
\circ t^{1/d}p_1)^{-1}\\
&\qquad\qquad\times (\Exp_G\circ t^{-1/d}(\Exp^\onemodd)^{[-1]}
\circ t^{1/d}p_1).
\end{split}
\ees
\eth
\noindent
In the $d=2$ special case, a better notation for the 
right-hand side is
\bes
\begin{split}
&-t^{1/2}\Tanh_G\circ\Arcsinh\circ t^{1/2}p_1\\
&\qquad+(\Sech_G\circ\Arcsinh\circ t^{1/2}p_1)(\Exp_G\circ 
t^{-1/2}\Arcsinh\circ t^{1/2}p_1).
\end{split}
\ees
The non-equivariant version of this special case is:
\beq
\begin{split}
1+\sum_{n\geq 1}
&\sum_{i\in\Z}\dim\WH_i(Q_n^{1\,\mathrm{mod}\,2}(G))\,(-t)^i
\frac{x^n}{|G|^n n!}\\
&=-t^{1/2}\tanh(\frac{1}{|G|}\arcsinh(t^{1/2}x))\\
&\qquad+\sech(\frac{1}{|G|}\arcsinh(t^{1/2}x))\exp(\frac{t^{-1/2}}{|G|}
\arcsinh(t^{1/2}x)).
\end{split}
\eeq
For any $n\geq 2$, let $\Pi_{B_n}^{(2)}$ denote the poset of (conjugate)
parabolic subsystems of a root system of type $B_n$ all of whose components
have even rank. As noted after Definition \ref{1modddef}, this is 
$Q_n^{1\,\mathrm{mod}\,2}(\{\pm 1\})$ if $n$ is even and
$Q_n^{1\,\mathrm{mod}\,2}(\{\pm 1\})\setminus\{\hat{1}\}$ if $n$ is odd.
So in calculating the Whitney homology of $\Pi_{B_n}^{(2)}$, the first
term on the right-hand side of \eqref{fourthwhitneyeqn} drops out, and
Theorem \ref{mainthm} becomes:
\beq \label{mainbneqn}
\begin{split}
1&+\bbone_{\{\pm 1\}}^1+\sum_{n\geq 2}\sum_{i\in\Z}
\ch_{\{\pm 1\}\wr S_n}(\WH_i(\Pi_{B_n}^{(2)}))\,(-t)^i\\
&=(\Sech_{\{\pm 1\}}\circ\Arcsinh\circ t^{1/2}p_1)(\Exp_{\{\pm 1\}}\circ 
t^{-1/2}\Arcsinh\circ t^{1/2}p_1).
\end{split}
\eeq
The non-equivariant version is:
\beq \label{noneqmainbneqn}
\begin{split}
1+\frac{x}{2}&+\sum_{n\geq 2}\sum_{i\in\Z}
\dim\WH_i(\Pi_{B_n}^{(2)})\,(-t)^i\frac{x^n}{2^n n!}\\
&=\sech(\frac{1}{2}\arcsinh(t^{1/2}x))\exp(\frac{t^{-1/2}}{2}
\arcsinh(t^{1/2}x)).
\end{split}
\eeq
By unpublished work of Rains, the polynomial
$\sum_{i}\dim\WH_i(\Pi_{B_n}^{(2)})\,(-t)^i$ equals the
Poincar\'e polynomial of the manifold $\overline{\mathcal{M}}_{B_n}(\R)$,
consisting of the real points of the
De Concini-Procesi compactification of the complex hyperplane
complement of type $B_n$. Hence \eqref{noneqmainbneqn} also
gives the generating function for these Poincar\'e polynomials.
It is almost but not quite true that $\WH_i(\Pi_{B_n}^{(2)})
\cong H^i(\overline{\mathcal{M}}_{B_n}(\R);\Q)$ as representations of the
wreath product (i.e.\ the Coxeter group $W(B_n)$); there is some twisting
analogous to the tensoring by the sign character in 
the type $A$ result \cite[Theorem 3.5]{rains}. Similar remarks apply
in type $D$, so it is worth recording the analogue of \eqref{mainbneqn}:
\beq \label{maindneqn}
\begin{split}
&\qquad 1+\bbone_{\{\pm 1\}}^1+\bbone_{\{\pm 1\}}^2+\sum_{n\geq 3}\sum_{i\in\Z}
\ch_{\{\pm 1\}\wr S_n}(\WH_i(\Pi_{D_n}^{(2)}))\,(-t)^i\\
&=(1+t\bbone_{\{\pm 1\}}^2)(\Sech_{\{\pm 1\}}\circ\Arcsinh\circ t^{1/2}p_1)
(\Exp_{\{\pm 1\}}\circ t^{-1/2}\Arcsinh\circ t^{1/2}p_1),
\end{split}
\eeq
and of \eqref{noneqmainbneqn}:
\beq \label{noneqmaindneqn}
\begin{split}
1&+\frac{x}{2}+\frac{x^2}{8}+\sum_{n\geq 3}\sum_{i\in\Z}
\dim\WH_i(\Pi_{D_n}^{(2)})\,(-t)^i\frac{x^n}{2^n n!}\\
&=(1+\frac{tx^2}{8})\,\sech(\frac{1}{2}\arcsinh(t^{1/2}x))
\exp(\frac{t^{-1/2}}{2}\arcsinh(t^{1/2}x)).
\end{split}
\eeq
These can be proved by the same method.

Finally, we deduce the following result from \eqref{suspendedeqn} and some
mild algebraic manipulation.
\bth
In $\Aa_G\otimes_\Q \Q[t,t^{-1}]$ we have the equation
\bes
\begin{split} 
1+\sum_{n\geq 1}
&\sum_{i\in\Z}\ch_{G\wr S_n}(\WH_i(Q_n^\zeromodd(G)))\,(-t)^i\\
&=\Exp_G+t
-(\sum_{j=0}^{d-1}t^{\frac{d-j}{d}}\Exp_G^{j\,\mathrm{mod}\,d}
\circ t^{1/d}p_1)\\
&\qquad\qquad\times(\Exp_G\circ (t^{-1}-1)L\circ
(\Exp^\zeromodd-1)\circ
t^{1/d}p_1).
\end{split}
\ees
\eth
\noindent
The non-equivariant version of the $d=2$ case is:
\beq
\begin{split}
1+\sum_{n\geq 1}
&\sum_{i\in\Z}\dim\WH_i(Q_n^{0\,\mathrm{mod}\,2}(G))\,(-t)^i\frac{x^n}{|G|^n n!}
=\exp(\frac{x}{|G|})+t\\
&-(t\cosh(\frac{t^{1/2}x}{|G|})+t^{1/2}\sinh(\frac{t^{1/2}x}{|G|}))
\cosh(t^{1/2}x)^{\frac{1}{|G|}(t^{-1}-1)}.
\end{split}
\eeq

\end{document}